\documentclass[12pt,reqno]{amsart}

\usepackage{AKstyle}

\begin{document}

\author{Jean Bourgain}
\thanks{Bourgain is partially supported by NSF grant DMS-0808042.}
\email{bourgain@ias.edu}
\address{IAS, Princeton, NJ}
\author{Alex Kontorovich}
\thanks{Kontorovich is partially supported by  NSF grants DMS-0802998 and DMS-0635607, and the Ellentuck Fund at IAS}
\email{avk@math.ias.edu}
\address{IAS and Brown University, Princeton, NJ}

\title[On a Theorem of Friedlander and Iwaniec]{
On a Theorem of Friedlander and Iwaniec
}

\date{\today}
\maketitle

\begin{abstract}
In \cite{FriedlanderIwaniec2009}, Friedlander and Iwaniec studied the so-called Hyperbolic Prime Number Theorem, which asks for an infinitude of elements $\g=\bigl( \begin{smallmatrix}
a&b\\ c&d
\end{smallmatrix} \bigr)
\in\SL(2,\Z)$ such that the norm squared
$$
\|\g\|^{2}=a^{2}+b^{2}+c^{2}+d^{2}=p,
$$ 
a prime. Under the Elliott-Halberstam conjecture, 
they 
proved the existence of such, as well as a  formula for their count, off by a constant from the conjectured asymptotic. In this note, we study 
the analogous question 
replacing the integers with the Gaussian integers. 
We prove unconditionally that 
for every odd $n\ge3$,
there is a 
$\g\in\SL(2,\Z[i])$ such that $\|\g\|^{2}=n$.
In particular, every prime 
is represented. 
 The proof is an application of Siegel's mass formula.
\end{abstract}

\section{Introduction}

The Affine Linear Sieve, introduced by Bourgain, Gamburd and Sarnak \cite{BourgainGamburdSarnak2006}, aims to produce prime points for functions on orbits of groups 
of morphisms of affine 
space. 
Friedlander and Iwaniec \cite{FriedlanderIwaniec2009} considered the case of the full modular group $\G=\SL(2,\Z)$, with the function being the norm-square.
Let $S$ be the set of norm-squares in $\G$, that is, 
$$
S:=\bigg\{n\in\Z_{+}\ : \ n=\|\g\|^{2}\text{ for some }\g\in\SL(2,\Z)\bigg\}.
$$
They proved, 
assuming
an approximation to
 the Elliott-Halberstam conjecture, that $S$ contains infinitely many  primes.\footnote{Moreover they gave a  formula for the count of norm-squares (with multiplicities), off by a constant from the conjectured asymptotic.}


Unconditionally, one can easily show the existence of $2$-almost primes in $S$. Indeed, for any $x\in\Z$, the parabolic elements 
$$
n_{x}:=\mattwo1x01
$$ 
are in $\G$, and their norm-square is $\|n_{x}\|^{2}=x^{2}+2$. Then Iwaniec's theorem \cite{Iwaniec1978} produces infinitely many $2$-almost primes in $S$.

In this note, we ask an analogous 
question, replacing the integers in $\SL(2,\Z)$ by the Gaussian integers, $\G=\SL(2,\Z[i])$. We prove unconditionally the following

\begin{thm}\label{thm:main}
The set
$$
S:=\bigg\{n\in\Z_{+}\ :\ n=\|\g\|^{2}\text{ for some }\g\in\SL(2,\Z[i])\bigg\}
$$
contains 
all odd integers $n\ge3$. 
In particular, it contains all primes.
\end{thm}

The proof, given in the next section, is an application of Siegel's mass formula \cite{Siegel1935}.
The argument 
is sufficiently delicate that it cannot
 replace the Gaussian integers above by the ring of integers of  another number field, even an imaginary quadratic extension (as suggested to us by John Friedlander), 
 see Remark \ref{rmk:1}.

%
%
We 
thank Zeev Rudnick and Peter Sarnak for conversations regarding this work.

\section{Sketch of the Proof}

For odd $n\ge3$ and   
$\g=\bigl(
\begin{smallmatrix}
a&b \\
c &d
\end{smallmatrix}
\bigr)
$
with $a=a_{1}+ia_{2}$, etc., the conditions $n=\|\g\|^{2}$ and $\g\in\SL(2,\Z[i])$ imply
\be\label{eq:sys1}
\threecase{
\|\g\|^{2}=
a_{1}^{2}+b_{1}^{2}+c_{1}^{2}+d_{1}^{2}+a_{2}^{2}+b_{2}^{2}+c_{2}^{2}+d_{2}^{2}=n
,
}{}
{
\Re(\det \g)=
a_{1} d_{1} -b_{1} c_{1} + b_{2} c_{2} - a_{2} d_{2}
=
1
}{}
{
\Im(\det \g)
=
a_{1} d_{2} + a_{2} d_{1}   - b_{1} c_{2} -b_{2} c_{1}
=
0
.
}{}
\ee
Changing variables 
\beann
a_{1}&\to& (y_{1}+y_{4})/2,\qquad\qquad
b_{1}\to (y_{3}+y_{2})/2,\\
c_{1}&\to& (y_{3}-y_{2})/2,\qquad\qquad
d_{1}\to (y_{1}-y_{4})/2,\\
a_{2}&\to& (y_{5}+y_{8})/2,\qquad\qquad
b_{2}\to (y_{7}+y_{6})/2,\\
c_{2}&\to& (y_{7}-y_{6})/2,\qquad\qquad
d_{2}\to (y_{5}-y_{8})/2,\\
\eeann
the system \eqref{eq:sys1} becomes
\be\label{eq:eqs}
\threecase{
 y_{3}^2 + y_{4}^2+ y_{5}^2 + y_{6}^2 =n-2,
}{}
{
y_{1}^2 + y_{2}^2 + y_{7}^2 + y_{8}^2 = n + 2,
}{}
{
y_{1} y_{5} + y_{2} y_{6} - y_{3} y_{7} - y_{4} y_{8}=0.
}{}
\ee
Write
$$
F=\bp1&&&\\&1&&\\&&1&\\&&&1\ep, \qquad\qquad G_{n}=\bp n+2 &\\&n-2\ep,
$$
and
$$
X=\bp 
y_{1} & y_{2} &-y_{7}& -y_{8}
\\
y_{5} &  y_{6} & y _{3} & y_{4}
\ep
,
$$
so that \eqref{eq:eqs} becomes
\be\label{eq:eqsMat}
X F \, ^{t}X=G_{n}.
\ee
\

Recall Siegel's \cite{Siegel1935} mass formula, cf. \cite[Appendix B, equations (3.10) to (3.17)]{Cassels1978}.
Clearly $F$ is 
positive definite and alone in its genus,  
and hence
the number $\cN(F,G_{n})$ of solutions $X$ to \eqref{eq:eqsMat}  is 
given by
\be\label{eq:cNDef}
\cN(F,G_{n})=\prod_{p\le\infty}\ga_{p}(F,G_{n})
,
\ee
where 
the 
local densities
$\ga_{p}$ are given as follows. For $p<\infty$, they are defined by
\be\label{eq:apDef}
\ga_{p}(F,G_{n}) = p^{-5t}\cdot \#\big\{X(\mod p^{t}): XF\,^{t}X\equiv G_{n}(p^{t})\big\}, 
\ee

\noindent 
for $t$ sufficiently large. For $p=\infty$, we have
$$
\ga_{\infty}(F,G_{n}) = 2\pi^{3} (n^{2}-4)^{1/2} 
.
$$
\begin{rmk}\label{rmk:1}
In complete generality, it is 
notoriously
difficult 
to compute the local densities $\ga_{p}$ 
and 
extract information such as non-vanishing,
see e.g. the formulae in \cite{Yang1998, Yang2004}.  
The main problem  being how large is ``sufficiently large'' for $t$ in \eqref{eq:apDef} with a given $p$.
In our special case of $\G=\SL(2,\Z[i])$, the literature is sufficient  to carry out the task.
\end{rmk}

For $p\neq2$, both the ramified and unramified local densities can be evaluated as in e.g. \cite[Theorem 2]{Kitaoka1983}. 
We turn first to the case $p$ is unramified,  $p\nmid (n^{2}-4)$. Then
$$
\ga_{p}(F,G_{n}) = \left(1-\frac1{p^{2}}\right)\left(1+{\chi_{p}\hskip-.04in
\left(
{4-n^{2}}
\right)
\over
p
}
\right)
,
$$
where $\chi_{p}
=\left(\frac \cdot p\right)$ 
is the quadratic character 
mod $p$.
%
For ramified primes $p\ge3$,
write $n+2=mp^{a}$ and $n-2=kp^{b}$ with $(mk,p)=1$. Assume $0\le a\le b$ (otherwise reverse their roles). Then if $a+b\equiv0(\mod 2)$,
$$
\ga_{p}(F,G_{n}) =
\frac{(p+1) \left(\left(p^{a+1}-1\right)
   \big(\chi_{p}\hskip-.03in\left({-mk
   }\right)-1\big)+(a+1)
   \left(p^2-1\right)  p^{(a+b)/2} \right)}{
    p^{3+(a+b)/2}}
   .
$$
Otherwise, if $a+b\equiv1(\mod 2)$, then
$$
\ga_{p}(F,G_{n}) =
\frac{
 (p+1)^2 \left(
(a+1)(p-1) p^{
   (a+b+1)/2}
    -
{\left(
p^{a+1}-1\right) 
}
   \right)}{
   p^{
   3+(a+b+1)/2}
   }
   .
$$
Inspection shows that these terms never vanish.

\

It remains to evaluate the dyadic density, $\ga_{2}$
. As shown by Siegel, see \cite{Reiner1956}, for $n$ odd (and hence $n^{2}-4$ odd), it is sufficient to evaluate \eqref{eq:apDef} for $t=3$, that is, compute the number of solutions $\mod 2^{3}=8$. 
One can compute explicitly that for any odd $n$, the number of solutions to \eqref{eq:apDef} $\mod 8$ is $49\,152$. Since $8^{5}=32\,768$, we have
$$
\ga_{2}(F,G_{n})=\frac32.
$$

In conclusion, the $\ga_{p}$'s never vanish so there are no local obstructions for odd $n$ to be represented, and hence the set $S$ of norm-squares in $\SL(2,\Z[i])$ contains all the  primes. (The prime $2$ is in $S$ since it is the norm squared of the identity matrix.)

\bibliographystyle{alpha}
\bibliography{../../AKbibliog}
\end{document}